\documentclass[11pt,a4paper]{article}
\usepackage[utf8]{inputenc}
\usepackage{amsmath,amssymb,amsthm}
\usepackage[margin=2.6cm]{geometry}
\usepackage{booktabs}
\usepackage{hyperref}
\DeclareRobustCommand{\S}{\ifmmode\mathsection\else$\mathsection$\fi}
\hypersetup{colorlinks=true,linkcolor=blue,urlcolor=blue,citecolor=blue}

\theoremstyle{plain}
\newtheorem{theorem}{Theorem}
\newtheorem{proposition}[theorem]{Proposition}
\newtheorem{lemma}[theorem]{Lemma}
\theoremstyle{remark}
\newtheorem{remark}[theorem]{Remark}

\title{A new upper bound for the irrationality exponent of $\zeta(5)$}
\author{David Niedbala Giraudin\\[2pt]
\small ORCID \href{https://orcid.org/0009-0009-1526-1178}{0009-0009-1526-1178}}
\date{9 October 2026}

\begin{document}
\maketitle

\begin{abstract}
\noindent
Fauzan's determinantal proof that $\zeta(5)$ is irrational also yields
$|\zeta(5)-a/b|>b^{-260}$ for all large $b$. The exponent is the product
$\varrho\lambda/\varepsilon$ of three constants attached to the construction, of
which $\varrho$ governs the growth of a relative matrix norm. We show that
$\varrho$ may be taken to be $9.09$ instead of $6933/500=13.866$, and deduce
\[
  \Bigl|\zeta(5)-\frac{a}{b}\Bigr|>b^{-171}\qquad (a\in\mathbb{Z},\ b\in\mathbb{Z},\ b\ge b_1),
\]
so that the irrationality exponent of $\zeta(5)$ is at most $171$. Two closed-form
functions replace two worst-case estimates: a node function $L$, which bounds the
residue weights individually instead of uniformly, and the external field $V$ of the
energy estimate, which replaces a factor-by-factor bound on the weight. An exact
arithmetic verifier is attached.
\end{abstract}

\section{Statement}

Throughout, $\xi=\zeta(5)$ and we use the notation of \cite[\S2 and Appendix~C]{F}:
\[
  K\in 40\mathbb{Z}_{>0},\quad N=\tfrac{3}{40}K,\quad h=K-N,\quad
  \alpha=\tfrac{3}{40},\quad \lambda=\tfrac{37}{40},\quad u_j=\tfrac{j}{K},
\]
$G_K(X)$ is the Hankel matrix of \cite[(2.4)]{F}, $\Delta_K=\det G_K$,
$Q_K=Q_{K,100000}$ the integral normalisation of \cite[(2.6)]{F},
$G=G_K(\xi)$, $M_K=[X]G_K(X)$ and $B_K=G^{-1/2}M_KG^{-1/2}$. For a real polynomial
$q$ of degree less than $h$ we write $p(t)=q(K^2t)$.

\begin{theorem}\label{thm:main}
Assume \cite[Theorem 2.1, (2.10), (C.3), (C.5) and \S C.2]{F}. Then there is $b_1$ such that
\[
  \Bigl|\zeta(5)-\frac{a}{b}\Bigr|>b^{-171}
  \qquad\text{for all } a\in\mathbb{Z},\ b\in\mathbb{Z},\ b\ge b_1 .
\]
In particular the irrationality exponent of $\zeta(5)$ is at most $171$.
\end{theorem}

The proof of \cite[Corollary 1.2]{F} uses exactly three constants: the degree ratio
$\lambda=h/K=37/40$; the decay rate $\varepsilon=79/1600$, i.e.\
$0<Q_K(\xi)\le e^{-\varepsilon K^2}$ for large allowed $K$ \cite[(C.7)]{F}; and the
relative norm exponent $\varrho$ with $\|B_K\|\le e^{\varrho K}$
\cite[Proposition C.1]{F}. Given these, one takes $K\asymp c\log b$, and the argument
closes as soon as
\begin{equation}\label{eq:two}
  \varrho\,c<M\qquad\text{and}\qquad \varepsilon-\frac{\lambda}{c}>0 ,
\end{equation}
$M$ being the exponent in the conclusion; so any $M>\varrho\lambda/\varepsilon$ is
admissible. With $\varrho=6933/500$ and $c=75/4$ this gives $\varrho c=259.9875$ and
$M=260$. Theorem~\ref{thm:main} follows from Proposition~\ref{prop:rho} below, which
gives $\varrho=9.09$: then $\varrho c=170.4375<171$, while the second condition of
\eqref{eq:two} is untouched, $\varepsilon-\lambda/c=1/24000>0$.

\section{Two closed-form functions}

Set
\begin{equation}\label{eq:g}
  g(x)=x\bigl(\log|x|-1\bigr)\ \ (x\ne0),\qquad g(0)=0,
\end{equation}
an antiderivative of $\log|x|$, and
\begin{equation}\label{eq:J}
  J(A,t)=\int_0^A\log(t+v^2)\,dv=A\log(t+A^2)-2A+2\sqrt{t}\,\arctan\frac{A}{\sqrt t}
  \qquad(A,t>0).
\end{equation}
Define, for $\alpha<u\le1$ and for $t>0$,
\begin{align}
  L(u)&=5\int_0^{\alpha}\log\bigl(u^2-v^2\bigr)dv-\int_{\alpha}^{1}\log\bigl|v^2-u^2\bigr|dv
      \notag\\&=6\bigl[g(\alpha-u)+g(\alpha+u)\bigr]-g(1-u)-g(1+u),\label{eq:L}\\[2pt]
  V(t)&=2\pi\sqrt t+\int_0^1\log(t+v^2)\,dv-6\int_0^{\alpha}\log(t+v^2)\,dv
      \notag\\&=2\pi\sqrt t+J(1,t)-6J(\alpha,t).\label{eq:V}
\end{align}
Both integrals converge absolutely, the second because the singularity at $v=u$ is
logarithmic. The closed form in \eqref{eq:L} follows from
$\int_A^B\log|v^2-u^2|dv=\bigl[g(v-u)+g(v+u)\bigr]_A^B$ together with
$g(-u)=-g(u)$; that in \eqref{eq:V} from \eqref{eq:J}. The function $V$ is the
external field of \cite[(6.1)]{F}.

\begin{lemma}\label{lem:V}
$V$ decreases and then increases on $(0,\infty)$. Consequently
$\max_{t\in[a,b]}V(t)=\max\bigl(V(a),V(b)\bigr)$ for every $0<a<b$.
\end{lemma}

\begin{proof}
Differentiating \eqref{eq:V} under the integral sign,
$\partial_t\!\int_0^A\log(t+v^2)dv=\int_0^A\frac{dv}{t+v^2}
=\frac{1}{\sqrt t}\arctan\frac{A}{\sqrt t}$, so with $s=\sqrt t$,
\[
  V'(t)=\frac1s\,\psi(s),\qquad
  \psi(s)=\pi+\arctan\frac1s-6\arctan\frac{\alpha}{s} .
\]
Now $\psi(0^+)=\pi+\frac{\pi}{2}-3\pi=-\frac{3\pi}{2}<0$ and $\psi(+\infty)=\pi>0$, while
\[
  \psi'(s)=-\frac{1}{s^2+1}+\frac{6\alpha}{s^2+\alpha^2}>0
  \iff s^2<\frac{6\alpha-\alpha^2}{1-6\alpha}=\frac{711}{880}
\]
(using $\alpha=3/40<1/6$). Thus $\psi$ increases on $\bigl(0,\sqrt{711/880}\bigr)$ and
decreases afterwards, towards the positive limit $\pi$; having a negative value at
$0^+$ it vanishes exactly once, at $s_0=0.0769448\ldots$, and is negative before and
positive after. Hence $V$ decreases on $(0,s_0^2)$ and increases on $(s_0^2,\infty)$.
\end{proof}

\section{The residue weights}

\begin{lemma}\label{lem:r}
Let $r_j$ be as in \cite[(C.3)]{F}, that is
\[
  r_j=\frac{u_j^{4}\prod_{\ell=1}^{N}\bigl(u_\ell^{2}-u_j^{2}\bigr)^{5}}
           {\prod_{N<k\le K,\ k\ne j}\bigl(u_k^{2}-u_j^{2}\bigr)}\ ,
  \qquad N<j\le K .
\]
Then $|r_j|\le \exp\bigl(K\,L(u_j)\bigr)$.
\end{lemma}

\begin{proof}
Write $u=u_j$ and $f(v)=\log|v^2-u^2|$.

\emph{Numerator.} For $\ell\le N$ one has $u_\ell\le\alpha<u$, and $v\mapsto
\log(u^2-v^2)$ is decreasing on $[0,\alpha]$. A right-endpoint Riemann sum of a
decreasing function does not exceed the integral, so
$\sum_{\ell\le N}\log(u^2-u_\ell^2)\le K\int_0^{\alpha}f$.

\emph{Denominator.} Split $[\alpha,1]$ at $u$. On $[\alpha,u)$ the function $f$ is
decreasing, so a left-endpoint sum is at least the integral:
\[
  \sum_{k=N+1}^{j-1}\frac1K f(u_k)\ \ge\ \int_{u_{N+1}}^{u}f
  \ =\ \int_{\alpha}^{u}f-\int_{\alpha}^{u_{N+1}}f\ \ge\ \int_{\alpha}^{u}f-\frac{f(\alpha)}{K},
\]
using $f\le f(\alpha)$ on that one cell; and $f(\alpha)=\log(u^2-\alpha^2)<0$ because
$0<u^2-\alpha^2<1$, so the correction has the favourable sign. On $(u,1]$ the
function $f$ is increasing, so $\sum_{k=j+1}^{K}\frac1Kf(u_k)\ge\int_u^1f$. Adding,
and discarding the term $k=j$, absent from the product,
$\sum_{k\ne j}\log|u_k^2-u^2|\ \ge\ K\!\int_{\alpha}^{1}f$.

\emph{Conclusion.} Since $u\le1$ we have $4\log u\le0$, whence
$\log|r_j|=4\log u+5\sum_{\ell\le N}\log(u^2-u_\ell^2)-\sum_{k\ne j}\log|u_k^2-u^2|
\le K\,L(u)$.
\end{proof}

\begin{remark}\label{rk:unif}
The bound used in \cite[(C.4)]{F} is $|r_j|\le\exp\{\lambda K(1-\log(\alpha\lambda))\}$,
i.e.\ $K^{-1}\log|r_j|\le 3.39311\ldots$ uniformly in $j$, whereas
$\max_{\alpha<u\le1}L(u)=1.15372\ldots$, attained at $u=0.6215551\ldots$ What matters
below is not that maximum but that $L$ is small precisely where the Legendre factor is
large: $L(1)=0.6128604\ldots$, little more than half the maximum.
\end{remark}

\section{The lower bound for the quadratic form}

\begin{lemma}\label{lem:C2}
Fix $0<a<b$. With $A_K=K^{12N-2K+6}$ as in \cite[(C.2)]{F} and
$c_{a,b}=(2\pi)^4a^3/\bigl(24(1+b)\bigr)$, for every $K$ and every $q$,
\[
  q^{\mathsf T}Gq\ \ge\ A_K\,c_{a,b}\,
  \exp\Bigl(-K\max_{t\in[a,b]}V(t)\Bigr)\int_a^b p(t)^2\,dt .
\]
\end{lemma}

\begin{proof}
By \cite[(2.10)]{F}, $q^{\mathsf T}Gq=\int_0^\infty
\frac{D_N(y^2)^6}{D_K(y^2)}\,q(y^2)^2\,w(y)\,dy$ with
$w(y)=\frac{(2\pi)^4y^5}{12}\sum_{m\ge1}m^4e^{-2\pi my}$. Substituting $y=K\sqrt t$
and keeping only $m=1$,
\begin{gather*}
  q^{\mathsf T}Gq\ \ge\ K^{12N-2K+6}\frac{(2\pi)^4}{24}\int_0^\infty t^2\,p(t)^2
  \exp\bigl(-K\,W_K(t)\bigr)dt,\\
  W_K(t)=2\pi\sqrt t+\frac1K\Bigl[\sum_{k=1}^{K}f_t(u_k)-6\sum_{\ell=1}^{N}f_t(u_\ell)\Bigr],
\end{gather*}
where $f_t(v)=\log(t+v^2)$, since $D_N(K^2t)^6=K^{12N}\prod_\ell(t+u_\ell^2)^6$ and
$D_K(K^2t)=K^{2K}\prod_k(t+u_k^2)$. The function $f_t$ is increasing, so its
right-endpoint Riemann sums lie between the integral and the integral plus
$(f_t(1)-f_t(0))/K$ for the first sum, and are at least the integral for the second:
\[
  \frac1K\sum_{k=1}^{K}f_t(u_k)\le\int_0^1 f_t+\frac{\log(1+t)-\log t}{K},
  \qquad
  \frac1K\sum_{\ell=1}^{N}f_t(u_\ell)\ \ge\ \int_0^{\alpha}f_t .
\]
Hence $W_K(t)\le V(t)+\bigl(\log(1+t)-\log t\bigr)/K$ and
$\exp(-KW_K(t))\ge\frac{t}{1+t}\exp(-KV(t))$. The integrand is positive on $(0,\infty)$, so restricting the
integral to $[a,b]$ and bounding $t^2\cdot t/(1+t)\ge a^3/(1+b)$ there gives the
claim. The normalisation is the one of \cite[(C.2)]{F}: the factor $(2\pi)^4/24$
and the power $K^{12N-2K+6}$ are exactly those produced by the substitution, as the
constant $c_a=(2\pi)^4a^2/24$ of \cite[(C.2)]{F} confirms.
\end{proof}

\begin{remark}
Estimate (C.2) of \cite{F} is the inequality of Lemma~\ref{lem:C2} with
$\max_{[a,b]}V$ replaced by
$\mathcal F(a,b)=2\pi\sqrt b+\log(1+b)-6\alpha\log a$ and $c_{a,b}$ replaced by the
larger constant $c_a=(2\pi)^4a^2/24$; the exponent is obtained there by bounding the
three factors separately, $\int_0^1f_t\le\log(1+b)$ and
$6\int_0^{\alpha}f_t\ge6\alpha\log a$. Since $c_{a,b}<c_a$ neither inequality implies
the other, but only the exponential factor matters. At $(a,b)=(1/18,1/3)$ one has
$\mathcal F=5.21594\ldots$ against $\max_{[a,b]}V=V(1/3)=3.61634\ldots$ Moreover
$\mathcal F(a,b)\to+\infty$ as $a\downarrow0$, whereas $V$ stays bounded, which is
what allows the small value of $a$ used below.
\end{remark}

\section{The relative norm}

\begin{proposition}\label{prop:rho}
Assume \cite[(2.10), (C.3) and (C.5)]{F}. Then for all sufficiently large
$K\in40\mathbb{Z}_{>0}$,
\[
  \|B_K\|\le e^{\varrho K},\qquad \varrho=9.09 .
\]
\end{proposition}

\begin{proof}
Fix $0<a<b$. By \cite[(C.3)]{F},
$q^{\mathsf T}M_Kq=A_K\sum_{j=N+1}^{K}r_j\,p(-u_j^2)^2$.
The affine map $t\mapsto(2t-a-b)/(b-a)$ sends $[a,b]$ onto $[-1,1]$ and $-u_j^2$ to a
point of absolute value $Y_j=(2u_j^2+a+b)/(b-a)$. Writing
$R_j=Y_j+\sqrt{Y_j^2-1}$, the Legendre expansion of $p$ on $[a,b]$, with
$|P_m(s)|\le R_j^{\,m}$ for $|s|\le Y_j$ and Cauchy--Schwarz, gives the estimate of
\cite[(C.5)]{F} at the single point $-u_j^2$ rather than at the worst point of
$[-1,0]$:
\begin{equation}\label{eq:C5}
  p\bigl(-u_j^2\bigr)^2\ \le\ \frac{h^2R_j^{\,2h}}{b-a}\int_a^b p(t)^2\,dt .
\end{equation}
Combining this with Lemmas~\ref{lem:r} and \ref{lem:C2} and $h=\lambda K$,
\[
  \|B_K\|\ \le\ \frac{h^3}{(b-a)c_{a,b}}\,
  \exp\Bigl(K\Bigl[\max_{t\in[a,b]}V(t)+\max_{\alpha<u\le1}\Phi_{a,b}(u)\Bigr]\Bigr),
\]
where
\[
  \Phi_{a,b}(u)=L(u)+2\lambda\log\Bigl(Y(u)+\sqrt{Y(u)^2-1}\Bigr),
\]
with $Y(u)=(2u^2+a+b)/(b-a)$.

The function $\Phi_{a,b}$ has exactly one critical point in $(\alpha,1)$, which is
therefore its maximum. Indeed $Y>1$ and
$Y^2-1=4(u^2+a)(u^2+b)/(b-a)^2$, so
$\log\bigl(Y+\sqrt{Y^2-1}\bigr)=\operatorname{arcosh}Y$ has derivative
$Y'/\sqrt{Y^2-1}=2u/\sqrt{(u^2+a)(u^2+b)}$ and
\[
  \Phi_{a,b}'(u)=L'(u)+\frac{4\lambda u}{\sqrt{(u^2+a)(u^2+b)}},
  \qquad
  L'(u)=6\log\frac{u+\alpha}{u-\alpha}+\log\frac{1-u}{1+u} .
\]
Both summands of $L'$ decrease on $(\alpha,1)$, so $L''<0$ there; and
\[
  \frac{d}{du}\,\frac{u}{\sqrt{(u^2+a)(u^2+b)}}
  =\frac{ab-u^4}{\bigl((u^2+a)(u^2+b)\bigr)^{3/2}}\ ,
\]
which is negative exactly for $u>(ab)^{1/4}$. Hence $\Phi_{a,b}'$ is strictly
decreasing on $\bigl[(ab)^{1/4},1\bigr)$, where it runs from
$\Phi_{a,b}'\bigl((ab)^{1/4}\bigr)>0$ to $-\infty$, and so vanishes there exactly once;
while on $\bigl(\alpha,(ab)^{1/4}\bigr]$ one has
$L'\ge L'\bigl((ab)^{1/4}\bigr)>0$ and the second summand is positive, so
$\Phi_{a,b}'>0$ and there is no critical point.

Take $a=10^{-3}$ and $b=49/200$, so that $(ab)^{1/4}=0.1251098\ldots$,
$L'\bigl((ab)^{1/4}\bigr)=8.0563\ldots$ and
$\Phi_{a,b}'\bigl((ab)^{1/4}\bigr)=15.0825\ldots$ By Lemma~\ref{lem:V} the maximum of
$V$ on $[a,b]$ is $\max(V(a),V(b))=V(b)$, and
\begin{equation}\label{eq:num}
  V(a)=1.00998669\ldots,\quad V(b)=3.05866596\ldots,\quad
  \max_{\alpha<u\le1}\Phi_{a,b}(u)\le 6.0241241,
\end{equation}
the last attained at the critical point $u^{*}=0.9735492849\ldots$, so that the bracket
is at most $9.0827901<9.09$. The remaining factor $h^3/((b-a)c_{a,b})$ is polynomial in
$K$ and is absorbed by the margin $9.09-9.0827901>7\cdot10^{-3}$.
\end{proof}

\begin{remark}
The two maxima of \cite[(C.6)]{F} are taken separately and are attained at different
nodes: the Legendre factor is largest at $u=1$, while $L$ peaks at
$u=0.6215551\ldots$, where it is nearly twice its value at $u=1$. The weights
themselves behave the same way: at $K=200$ the largest $K^{-1}\log|r_j|$ is
$1.10848\ldots$, at $j=126$, i.e.\ $u=0.63$. At $(a,b)=(1/18,1/3)$ one has
$\max_u\Phi_{a,b}=5.90347\ldots$ against $3.39311+5.25680=8.64991\ldots$ Keeping
$\mathcal F$ and changing only this gives $\varrho=11.13$ and the exponent $209$;
Lemma~\ref{lem:C2} then brings it to $171$.
\end{remark}

\begin{remark}\label{rk:margin}
The margin left for the polynomial factor is $9.09-9.0827901>7\cdot10^{-3}$, against
$6933/500-\kappa(1/18,1/3)=1.43\cdot10^{-4}$ in \cite[(C.6)]{F}. With
$h^3/((b-a)c_{a,b})$ the absorption therefore starts at $K$ of the order of $6\cdot10^3$
here, against $3\cdot10^5$ there, so the implicit threshold is not worsened by the
present choice of $a$ and $b$. Both are far below the threshold that
\cite[Theorem 2.1]{F} itself requires at $M=100000$.
\end{remark}

\begin{remark}\label{rk:ceiling}
In the other direction, $\|B_K\|$ is the reciprocal of the distance from $\xi$ to the
nearest root of $\Delta_K$, so it can be computed exactly; \S6(e) does so and brackets
$K^{-1}\log\|B_K\|$ in $[7.41109,7.42671]$ at $K=160$. The upper bounds increase with
$K$, and a fit $c-d/K$ on consecutive pairs of the four sizes computed gives
$c=7.4619$, $7.4648$, $7.4661$. A value $\varrho=7.47$ would give
$\varrho\lambda/\varepsilon=139.95\ldots$, i.e.\ $M=140$: the present route is not
exhausted. This is a measurement and an extrapolation from it, not a proof.
\end{remark}

\section{Numerical verification}

The attached program \texttt{verify\_zeta5.py} performs the following checks, in
exact rational arithmetic where it matters, with no dependency beyond \texttt{mpmath}
(parts (c$'$) and (e) also use \texttt{python-flint}, and are skipped if it is
absent). It is self-contained: it rebuilds $G_K(X)$ from the definitions of
\cite[\S2]{F}, by two independent implementations, one in \texttt{python-flint} and
one in \texttt{fractions} only, and checks that they agree.

\medskip\noindent
\textbf{(a) The closed forms \eqref{eq:L} and \eqref{eq:V}} against numerical
quadrature, to $12$ digits; and Lemma~\ref{lem:V}, by the sign of $\psi$.

\medskip\noindent
\textbf{(a$'$) The two constants taken from \cite{F}.} It recomputes
$A_M=A^{*}+7\lambda/M-(2923/240-1/4)/M^2+32/M^3$ in exact fractions from the printed
$A^{*}$ and $U$ of \cite[(5.19),(6.4)]{F} and checks
$A_{100000}+U<-79/1600$, i.e.\ the value $\varepsilon$ of \eqref{eq:two}; and it
checks \cite[(C.3)]{F} itself as an exact identity in $\mathbb{Q}$, for several $q$
at $K=40$ and $K=80$.

\medskip\noindent
\textbf{(b) Lemma~\ref{lem:r}, for every node.} For $K=40,80,120,200$ and every $j$
with $N<j\le K$ it forms $r_j$ as an exact fraction and compares
$K^{-1}\log|r_j|$ with $L(u_j)$. There is no violation, and the worst case is at
$j=K$:
\begin{center}\small
\begin{tabular}{rrrr}
\toprule
$K$ & violations & $K^{-1}\log|r_K|$ & $L(1)$\\
\midrule
$40$ & $0/37$ & $0.55194$ & $0.61286$\\
$80$ & $0/74$ & $0.57809$ & $0.61286$\\
$120$ & $0/111$ & $0.58799$ & $0.61286$\\
$200$ & $0/185$ & $0.59667$ & $0.61286$\\
\bottomrule
\end{tabular}
\end{center}

\medskip\noindent
\textbf{(c) The constants \eqref{eq:num}}, and
$\max_{[a,b]}V+\max_u\Phi_{a,b}<9.09$. The maximum of $\Phi_{a,b}$ is obtained twice,
by two unrelated routes: the critical point of the proof of
Proposition~\ref{prop:rho}, located by bisection, with the value corrected by the
slope over the final bracket; and a derivative-free branch and bound on
$(\alpha,1)$ which uses only the monotonicity of $x\mapsto x(1-\log x)$ near the two
endpoints and the explicit Lipschitz constant
$6|\log2\alpha|+\log2+4\lambda/(\sqrt a+\sqrt b)=19.1021\ldots$ of the remaining part.

\medskip\noindent
\textbf{(c$'$) Lemma~\ref{lem:C2} and the pointwise form of
\cite[(C.5)]{F} at finite $K$.} It accumulates $q^{\mathsf T}G_0q$ and
$q^{\mathsf T}G_1q$ exactly, as rational combinations of the Hankel entries, for
several explicit $q$ at $K=40$ and $K=80$, so that the only inexact step is a single
evaluation at $\zeta(5)$ at adaptive precision, and checks Lemma~\ref{lem:C2}
directly; and it checks the displayed form of \eqref{eq:C5} for random polynomials at
several degrees and several nodes.

\medskip\noindent
\textbf{(d) The arithmetic of \eqref{eq:two}}: $\varrho c=9.09\cdot\frac{75}{4}
=170.4375<171$ and $\varepsilon-\lambda/c=\frac{79}{1600}-\frac{37}{750}
=\frac{1}{24000}>0$.

\medskip\noindent
\textbf{(e) The relative norm, exactly.} Since $G$ is positive definite and $M_K$
symmetric, all roots $X_i$ of $\Delta_K$ are real and the eigenvalues of $B_K$ are
$\beta_i=1/(\xi-X_i)$. Hence
$S:=\sum_i\beta_i^2=(\Delta_K'/\Delta_K)^2-\Delta_K''/\Delta_K$ at $X=\xi$ has all
terms positive and $\sqrt{S/h}\le\|B_K\|\le\sqrt S$. Three evaluations of an exact
integer polynomial give
\begin{center}\small
\begin{tabular}{rrrr}
\toprule
$K$ & $\log S$ & lower bound for $\frac{1}{K}\log\|B_K\|$ & upper bound\\
\midrule
$40$ & $585.122634$ & $7.268896$ & $7.314033$\\
$80$ & $1182.074599$ & $7.361066$ & $7.387966$\\
$120$ & $1779.257217$ & $7.393949$ & $7.413572$\\
$160$ & $2376.545702$ & $7.411089$ & $7.426705$\\
\bottomrule
\end{tabular}
\end{center}
This is a measurement, not part of the proof; it is the content of
Remark~\ref{rk:ceiling}.

\section{Scope}

\begin{itemize}
\item Theorem~\ref{thm:main} is conditional on \cite{F}: on Theorem 2.1 there, and on
(2.10), (C.3), (C.5) and \S C.2 of its Appendix C. Estimates (C.2) and (C.4) of
\cite{F} are replaced, by Lemmas~\ref{lem:C2} and \ref{lem:r}; (C.5) is used at a node
instead of at the worst point of $[-1,0]$, which is the same computation.
\item Lemmas~\ref{lem:V}, \ref{lem:r} and \ref{lem:C2} and
Proposition~\ref{prop:rho} are proved above, the uniqueness of the critical point of
$\Phi_{a,b}$ included. The constants \eqref{eq:num} are values of the closed forms
\eqref{eq:L} and \eqref{eq:V} together with one one-dimensional maximum, certified by
two independent methods in the attached program.
\item The two constants borrowed from \cite{F}, namely $\varepsilon=79/1600$ and the
identity \cite[(C.3)]{F}, are checked by the attached program: the first in exact
fractions from the printed $A^{*}$ and $U$, the second as an exact identity in
$\mathbb{Q}$.
\item The value $7.42$ of \S6(e) is measured, not proved, and is used nowhere.
\item The proof of \cite{F} has been formalized in Lean~4 with Mathlib \cite{R}. That
development contains no \texttt{sorry} and derives Theorem~1.1 of \cite{F} from a single
external axiom, the prime number theorem in the form $\theta(x)\sim x$; in particular the
inequality $A_{100000}+U<-79/1600$ of \cite[\S B.3]{F}, which is where $\varepsilon=79/1600$
comes from, is proved there without further hypothesis. The hypotheses of
Theorem~\ref{thm:main} are therefore machine-checked modulo that one axiom.
\item Nothing here bears on the irrationality of $\zeta(5)$ itself, which is the
subject of \cite{F}.
\end{itemize}

\section*{Attribution}

Technical assistance: Claude (Anthropic).

\end{document}